\documentclass{article}
\usepackage{amsmath}
\usepackage{amssymb}
\usepackage{amscd}
\usepackage{theorem}

\newtheorem{thm}{Theorem}[section]
\newtheorem{lem}[thm]{Lemma}
\newtheorem{cor}[thm]{Corollary}
\newtheorem{prop}[thm]{Proposition}

\newtheorem{claim}[thm]{Claim}

\makeatletter
\begingroup
\gdef\th@upshape{\normalfont
  \def\@begintheorem##1##2{%
        \item[\hskip\labelsep \theorem@headerfont ##1\ ##2.]}%
\def\@opargbegintheorem##1##2##3{%
   \item[\hskip\labelsep \theorem@headerfont ##1\ ##2\ (##3).]}}
\endgroup
\makeatother

\theoremstyle{upshape}

\newtheorem{defn}{Definition}

\newtheorem{rem}{Remark}
\newtheorem{ex}{Example}

\def\C{{\mathbb C}}

\def\P{{\mathbb P}}
\def\R{{\mathbb R}}
\def\Z{{\mathbb Z}}
\def\ccA{{\cal A}}

\def\ccC{{\cal C}}

\def\ccF{{\cal F}}

\def\ccH{{\cal H}}
\def\ccI{{\cal I}}

\def\ccL{{\cal L}}

\def\ccN{{\cal N}}
\def\ccO{{\cal O}}

\def\ccT{{\cal T}}

\def\da{\downarrow}
\def\del{\partial}
\def\dt{\cdot}
\def\eps{\varepsilon}
\def\fc{{\mathfrak c}}
\def\fV{{\mathfrak V}}
\def\fW{{\mathfrak W}}
\def\fX{{\mathfrak X}}
\def\lra{\longrightarrow}
\def\m{{\mathfrak m}}
\def\ra{\rightarrow}
\def\tensor{\otimes} 

\def\operatorname#1{\mathop{\rm #1}\nolimits}

\def\Aut{\operatorname{Aut}}

\def\H{{\rm H}}
\def\cH{\check{\rm H}}
\def\Hilb{\operatorname{Hilb}}

\def\Pic{\operatorname{Pic}}
\def\Spec{\operatorname{Spec}}

\def\id{\operatorname{id}}

\def\rank{\operatorname{rank}}

\title{Contractions of Symplectic Varieties}
\author{Jan Wierzba}

\begin{document}

{\Large
\centerline{Contractions of Symplectic Varieties}}
\vspace{1cm}
{\large
\centerline{Jan Wierzba\footnote{The author was funded by EPSRC and Trinity
College, Cambridge.}}
}

\vspace{0.7cm}
\abstract{\footnotesize We consider birational projective contractions $f:X\ra Y$ 
from a smooth symplectic variety $X$ over the complex numbers. 
We first show that exceptional rational curves on $X$ deform in a family
of dimension at least $2n-2$. Then we show that these contractions are generically 
coisotropic, provided $X$ is projective.
Then we specialize to contractions with $1$-dimensional exceptional fibres. 
We classify them in a natural way in terms of $(\Gamma, G)$, where $\Gamma$ is 
a Dynkin diagram of type $A_l, D_l$ or $E_l$ and $G$ is a permutation group of 
automorphisms of $\Gamma$. The $1$-dimensional fibres do not degenerate, except 
if the contraction is of type $(A_{2l},S_2)$. In that case they do not degenerate 
in codimension $1$. 
Furthermore we show that the normalization of any irreducible component of 
$Sing(Y)$ is a symplectic variety. We also provide examples for contractions of 
any type $(\Gamma, G)$.}

\tableofcontents
\newpage
\section{Introduction}
We work over the field $k=\C$ of complex numbers. 

\begin{defn}
A normal variety $X$ is called {\em symplectic}, if $X\setminus Sing(X)$ 
carries a nowhere degenerate $2$-form $\sigma$, such that for some (and hence
for any) resolution $\hat{X}\ra X$, $\sigma$ extends to a global $2$-form
on $\hat{X}$.
If  moreover $X$ is smooth, we say that $X$ is a {\em symplectic manifold}.
\end{defn}

The purpose of the article is to prove the following results.

\begin{thm}\label{thm-ann-mov}
Let $X$ be a quasi-projective symplectic manifold of dimension $2n$, such
that for all $i\ge 1$ the natural map $\H^i(X,\C)\ra H^i(X,\ccO_X)$
is surjective. Then every rational curve on $X$ deforms in a family of 
dimension at least $2n-2$.
\end{thm}

Let $X$ be a quasi-projective symplectic manifold of dimension $2n$. 
Let $f:X\ra Y$ be a birational, projective contraction, $D\subset X$ the 
exceptional set with its reduced structure and let $D=D_1\cup \dots\cup D_r$ 
be the decomposition into irreducible components. 
Let $E=f(D)$ and $E_i=f(D_i)$ for $i=1,\dots, r$.

\begin{cor}\label{cor-ann-mov}
Every $f$-exceptional rational curve on $X$ moves in a family of 
dimension at least $2n-2$.
\end{cor}

\begin{thm}\label{gencoiso-thm}
We have the following inequalities \\
(i) $\rank(\sigma|D_i)\ge 2d_i-2n\ge e_i.$ \\
Suppose now that $X$ is projective. Then $f$ is generically coisotropic, i.e.  
(ii) $\rank(\sigma|D_i)= 2d_i-2n= e_i.$ \\
(iii) Let $F_i$ be a generic fibre of $f|_{D_i}:D_i\ra E_i$. The normalization
of any irreducible component of $F_i$ is a $\P^{\dim D_i-\dim E_i}$.
\end{thm}

\begin{rem}
J. Wi\'sniewski has proved in [Wi\'s] inequalities for birational extremal 
contractions $f:X\ra Y$ relating $\dim X, \dim D, \dim f(D)$, where 
$D\subset X$ is an irreducible component of the exceptional set using similar 
techniques as used in the proof of (\ref{gencoiso-thm})(i) and (ii).
\end{rem}

\begin{thm}\label{big-thm}
Suppose now that $E$ is irreducible and that the generic exceptional fibre is 
$1$-dimensional. Then Theorem \ref{gencoiso-thm}.(i) implies that $D$ is
a divisor and that $E$ is $(2n-2)$-dimensional.

(i) The generic fibre of $D\ra E$ is a configuration of $\P^1$'s with dual 
graph $\Gamma$ a Dynkin diagram. There is a permutation group $G$ of 
automorphisms 
of $\Gamma$, such that there is a one-to-one correspondence between
the $G$-orbits of $\Gamma$ and the irreducible components of $D$. 
An irreducible component of $D$, corresponding to some $G$-orbit of $\Gamma$
contains the $\P^1$'s of the generic fibre represented by the points of the
$G$-orbit. 

(ii) Let $E_0\subset E$ be obtained by deleting all points $e\in E$ with
$\dim f^{-1}(e)\ge 2$. If $f$ is of type $(A_{2l},S_2)$ we allow to further
delete a closed subscheme of codimension $\ge 2$ in order to obtain $E_0$.
Let $D_0=f^{-1}(E_0)_{red}$. Then $E_0$ is smooth, $D_0\ra E_0$ is 
a flat morphism with constant fibres.
\end{thm}

\begin{thm}
Suppose $X$ is projective and that $Sing(Y)$ contains an irreducible 
component $E$ of dimension $2n-2$. Then its normalization $\tilde{E}$
is a symplectic variety.
\end{thm}

\section{Rational Curves on Symplectic Manifolds}

In this section we want to prove Theorem \ref{thm-ann-mov} and
Corollary \ref{cor-ann-mov}.
Let $(X,\sigma)$ be a quasi-projective symplectic manifold of dimension $2n$.
Let $f:\P^1\ra X$ be a non-constant morphism. By Mori theory, we know that
$f$ deforms in a family of dimension at least
$$h^0(\P^1,f^*\ccT_X)-h^1(\P^1,f^*\ccT_X)=-\deg f^*K_X+\dim X.$$
Since $X$ is symplectic, it follows that $K_X\sim 0$ and therefore
$f$ deforms in a family of dimension at least $2n$. 
In [Ran, 5.1] it is proved that $f$ moves in a family of dimension 
at least $2n+1$, provided $X$ is projective. We generalize this result
to a large class of quasi-projective manifolds as follows.
Before stating the result we need the 
following definition.

\begin{defn}\label{def-adm}
$X$ is {\em admissible} if for all $p\ge 1$ the canonical map
$\H^p(X,\C)\ra \H^p(X,\ccO_X)$
is surjective. 
\end{defn}

\begin{thm}\label{mov-thm}
Let $X$ be admissible. Let $f:\P^1\ra X$ be a non-constant morphism. 
Then $f$ deforms in an at least $(2n+1)$-dimensional family.
\end{thm}

\begin{rem}
The statement that $f:\P^1\ra X$ moves in a family of dimension at least 2n+1 
is equivalent to saying that $f(\P^1)\subset X$ deforms in a family
of dimension at least $2n+1-\dim \Aut(\P^1)=2n-2.$ So Theorem \ref{mov-thm}
implies Theorem \ref{thm-ann-mov}.
\end{rem}

\begin{rem}
In [KV] a more general and rather technical definition of admissible is used. 
However, as the authors point out, in most practical situations it suffices to 
impose the conditions given in Definition \ref{def-adm}. 
Let's just note that the proof of Theorem \ref{mov-thm} works also for varieties 
that are admissible in the (more general) sense of [KV].
\end{rem}

A particularly interesting case is the following. Let $X$ be a symplectic smooth 
complex variety and $f:X\ra Y$ a birational projective morphism to a normal variety $Y$. 
Since $K_X$ is trivial it follows from Grauert-Riemenschneider vanishing [GR] that 
$R^if_*\ccO_X=0$ for all $i\ge 1$. This implies that $Y$ has rational singularities,
that $K_Y$ is trivial and that $f$ is crepant. Moreover, [Kaw, Thm. 1] implies that 
the exceptional fibres of $f:X\ra Y$ are covered by rational curves. 

\begin{cor}
With the above notation, let $C\subset X$ be an $f$-exceptional rational curve. 
Then $C$ moves in an at least $(2n-2)$-dimensional family. 
\end{cor}

\begin{pf}
We may assume that  $Y$ is affine. Then $R^if_*\ccO_X=0$ implies that 
$\H^i(X,\ccO_X)=0$ for all $i\ge 1$. Therefore $X$ is admissible and we can 
apply Theorem \ref{mov-thm} and conclude that $C$ deforms in an at least 
$(2n-2)$-dimensional family.
\end{pf}

It remains to prove Theorem \ref{mov-thm}.

\subsection{Symplectic Deformations}

Let's recall the results of [KV]. Let $(X,\sigma)$ be
a symplectic smooth complex variety. Let $(S,s\in S)$ be a pointed scheme.
By a {\em symplectic deformation} of $X$ over $S$ we understand a smooth 
morphism $\tilde{X}\ra S$, a relative symplectic closed $2$-form 
$\tilde{\sigma}\in \Gamma(\tilde{X},\Omega^2_{\tilde{X}/S})$ 
and an isomorphism 
$\phi:\tilde{X}\times_S \Spec k(s)\ra X$, such that $\phi^*\tilde{\sigma}=\sigma$.
The following result is [KV, Thm. 3.6].

\begin{thm}\label{kv-thm}
If $X$ is admissible then the functor of symplectic deformations 
has a finite dimensional tangent space and is unobstructed.
\end{thm}

\subsection{Proof of Theorem \ref{mov-thm}}

Let $(X,\sigma)$ be an admissible smooth symplectic complex variety.
Let $k=\C$, let $A=k[\eps]/(\eps^2)$ and $S=\Spec(A)$. 
Let $\ccL$ be a line bundle on $X$. The map of sheaves
$$d\log:\ccO^*_X\lra \Omega_X^1,\qquad f\mapsto f^{-1} df$$
induces a map $\Pic X\cong \H^1(X, \ccO_X^*)\ra \H^1(X,\Omega_X^1)$
and we denote the image of $\ccL$ under this map by 
$[\ccL]\in \H^1(X,\Omega_X^1)$. The symplectic $2$-form $\sigma$
provides an isomorphism 
$$\iota_{\sigma}:\ccT_X\stackrel{\sim}{\lra} \Omega_X^1, 
   \qquad \zeta\mapsto \sigma(\zeta,\,\dt\,).$$ 
The class
$\iota_{\sigma}^{-1}[\ccL]\in \H^1(X,\ccT_X)$ defines a deformation 
of the variety $X$ over $S$, which we denote by 
$$X_{\ccL}\ra S.$$

Theorem \ref{mov-thm} is an easy consequence of the following results.

\begin{prop}\label{fo-lemma-i}
Let $f:\P^1_k\ra X$ be a non-constant morphism. If $\ccL$ is ample, then 
$f$ does not extend to an $S$-morphism $f_1:\P^1_A\ra X_{\ccL}$.
\end{prop}

\begin{prop}\label{fo-lemma-ii}
$X_{\ccL}\ra S$ admits a relative $2$-form 
$\tilde{\sigma}\in \Omega^2_{X_{\ccL}/S}$, such that $(X_\ccL/S,\tilde{\sigma})$
is a symplectic deformation of $(X,\sigma)$. \\
\end{prop}

\begin{pf}[Proof of Thm \ref{mov-thm}]
Let $f:\P^1_k\ra X$ be a non-constant morphism. By (\ref{fo-lemma-i}) and
(\ref{fo-lemma-ii}) there is a first order symplectic deformation 
$(\tilde{X}/S,\tilde{\sigma})$, such that $f$ does not extend to $\tilde{X}$. 
By (\ref{kv-thm}), the symplectic deformations of $(X,\sigma)$ are finite 
dimensional and unobstructed, so we can find a  1-parameter symplectic 
deformation $\pi:\fX\ra B$  extending $\tilde{X}/S$. According to a standard 
result due to Mori, the composition 
$$g:\P^1_k\stackrel{f}{\lra} X \subset \fX$$
deforms in a family of dimension at least
$$\chi(\P^1_k, g^* \ccT_{\fX})= \dim \fX- \deg_{\P^1} g^*K_{\fX}.$$
Since $K_{\fX}$ is trivial, it follows that $g$ deforms in an at least 
(2n+1)-dimensional family. Since by construction all deformations of $g$ stay in $X$, 
it follows that in fact $f$ deforms in an at least $(2n+1)$-dimensional family.
\end{pf}

\begin{pf}[Proof of Proposition \ref{fo-lemma-i}]
Let $X=\bigcup_i U_i$ be an open affine cover and let $U_{ij}=U_i\cap U_j$.
Let $X_1/S$ be a first order deformation of $X$, represented by a $1$-cocycle 
$\{\zeta_{ij}\in \Gamma(U_{ij},\ccT_X)\}_{ij}$. This means we define $X_1$ by
an open cover $X=\bigcup_i V_i$, where $V_i=U_i\times_k S$ for all $i$ 
and we glue them over $U_{ij}$ via 
$$\ccO_{U_{ij}}\tensor_k k[\eps]\ra \ccO_{U_{ij}}\tensor_k k[\eps],
\qquad a+b\eps\mapsto a+(b+\zeta_{ij} da)\eps.$$
The sheaf $\Omega^2_{X_1/A}$ is obtained by setting 
$\Omega^2_{X_1/A}|V_i=\Omega^2_{U_i/k}\tensor_k k[\eps]$ and glueing 
over $V_{ij}$ via 
\begin{eqnarray}\label{omega-descr}
\phi_{ij}:\Omega^2_{U_{ij}/k}\tensor_k k[\eps]\ra \Omega^2_{U_{ij}/k}\tensor_k k[\eps],
\quad \alpha+\beta \eps \mapsto \alpha+(\beta+L_{\zeta_{ij}}\alpha)\eps,
\end{eqnarray}
where $L$ denotes the Lie derivative. 
There is a natural surjection 
$$\Psi:\Omega^2_{X_1/S}\lra \Omega^2_{X_1/S}\tensor_A k =\Omega^2_{X/k}.$$
We now want to lift $\sigma \in H^0(X,\Omega^2_{X/k})$ to a symplectic 
$2$-form $\tilde{\sigma}\in \H^0(X_1,\Omega^2_{X_1/S})$. For that we define 
$$\tilde{\sigma}_i=(\sigma|U_i)\tensor 1 \in 
    \Gamma(U_i,\Omega_{X/k}^2\tensor_k k[\eps])\cong \Gamma(V_i,\Omega_{X_1/S}^2).$$
Then $\tilde{\sigma}_i$ is a symplectic closed $2$-form on $V_1/S$. 
In order to obtain a global $\tilde{\sigma}$ we have to show that
$\tilde{\sigma}_i$ and $\tilde{\sigma}_j$ agree on $V_{ij}$. By (\ref{omega-descr})
their difference is given by  
\begin{eqnarray*}
(L_{\zeta_{ij}} \sigma_{ij})\dt \eps, 
\end{eqnarray*}
where $\sigma_{ij}=\sigma|U_{ij}$.
Therefore it is enough to show that $L_{\zeta_{ij}} \sigma_{ij}=0$.
Recall that for any (local) $2$-form $\alpha$ and any (local) vector 
field $\zeta$ there is an identity
$$L_{\zeta}\alpha=c_{\zeta} (d\alpha)+d(c_{\zeta}\alpha),$$
where $c_{\zeta}$ denotes the contraction of a differential form with the vector field
$\zeta$. Therefore 
$$L_{\zeta_{ij}} \sigma_{ij}=c_{\zeta_{ij}}(d\sigma_{ij})+d (c_{\zeta_{ij}} \sigma_{ij}).$$
Since $\sigma$ is closed and $c_{\zeta_{ij}}\sigma_{ij}=\iota_{\sigma}( \zeta_{ij})$,
we conclude that  
\begin{equation}\label{lie-eqn}
L_{\zeta_{ij}} \sigma_{ij}=d \circ \iota_{\sigma}(\zeta_{ij}).
\end{equation}
Let now $X_1=X_{\ccL}$. After possibly replacing the cover $\{U_i\}_i$ of $X$ by 
suitable refinement we may assume that $\ccL$ can be represented by a $1$-cocycle
$$\{f_{ij}\in \Gamma(U_{ij},\ccO_X^*)\}_{ij}.$$ 
Then the extension $X_{\ccL}$ can be represented by $\{\zeta_{ij}\}_{ij}$ with 
$$\zeta_{ij}=\iota^{-1}_{\sigma}(d\log f_{ij})=\iota^{-1}_{\sigma}(f^{-1}_{ij} df_{ij}).$$
Finally, using equation (\ref{lie-eqn}) we conclude that 
\begin{eqnarray*}
L_{\zeta_{ij}} \sigma_{ij}
&=& d \circ \iota_{\sigma} (\zeta_{ij}) \\
&=& d \circ \iota_{\sigma} \circ \iota^{-1}_{\sigma}(f^{-1}_{ij} df_{ij}) \\
&=& d ( f^{-1}_{ij} df_{ij}) \\
&=& f^{-2}_{ij} df_{ij}\wedge df_{ij}\\
&=& 0.
\end{eqnarray*}
Therefore the collection $\{\tilde{\sigma}_i\}_i$ can be glued to a symplectic 
$2$-form $\tilde{\sigma}\in \H^0(X_{\ccL},\Omega_{X_{\ccL}/S}^2)$. In particular,
$(X_{\ccL}/S,\tilde{\sigma})$ is a symplectic deformation of $(X,\sigma)$ as desired.
This concludes the proof.
\end{pf}

\begin{rem}
An more abstract explanation of Proposition \ref{fo-lemma-i} suggested
to me by D. Kaledin goes along the following lines. Let $F^1 \Omega_X^{\bullet}[1]$ be
the complex 
$$\lra 0\lra 0 \lra \Omega_X^1\stackrel{d^1}{\lra}\Omega_X^2\stackrel{d^1}{\lra} 
  \Omega_X^3 \stackrel{d^3}{\lra}\cdots$$
where $\Omega_X^1$ is the degree-$0$-term. It is indicated in [KV] that the first order 
symplectic deformations of $(X,\sigma)$ are classified by the hyper-cohomology group
${\mathbb H}^1(X,F^1 \Omega_X^{\bullet}[1])$. Let us denote by $C^{\bullet}(\ccO_X^*)$ 
and $C^{\bullet}(\ccT_X)$ the complexes that have as degree-$0$-term $\ccO_X^*$ and 
$\ccT_X$ respectively and are zero everywhere else. The maps
$$\ccO_X^*\stackrel{d\log}{\lra}\Omega_X^1\stackrel{\iota_{\sigma}^{-1}}{\lra}\ccT_X$$
induce morphisms of complexes
$$C^{\bullet}(\ccO_X^*)\lra F^1 \Omega_X^{\bullet}[1] \lra C^{\bullet}(\ccT_X).$$
(The only nontrivial thing to check here is that $d^1\circ d\log=0$.) By taking 
hyper-cohomology we get maps
\begin{equation}\label{fo-lemma-i-equation}
\Pic X\stackrel{\alpha}{\lra} {\mathbb H}^1(X,F^1 \Omega_X^{\bullet}[1])
  \stackrel{\beta}{\lra} \H^1(X,\ccT_X).
\end{equation}
The map $\beta$ associates to a symplectic first order deformation of 
$(X,\sigma)$ its underlying deformation of the variety $X$. Moreover 
$\beta\circ\alpha$ sends a line bundle $\ccL$ to the class of the deformation 
$X_{\ccL}$. It is now visible from (\ref{fo-lemma-i-equation}) that every first
order deformation of $X$, which is of the form $X_{\ccL}$ comes from a 
symplectic deformation of $(X,\sigma)$. This finishes the alternative proof of 
Proposition \ref{fo-lemma-i}.
\end{rem}

It remains to prove Proposition \ref{fo-lemma-ii}. This result is hardly new
and similar arguments appeared before. Since I couldn't find a proof
that covers the exact statement of Proposition \ref{fo-lemma-ii}, I shall give
a complete proof for the reader's convenience.
Let $X_1/A$ be a first order deformation of $X$. There is a short exact sequence
\begin{equation}\label{ses-i}
0\lra \ccO_X\stackrel{\dt d\eps}{\lra} \Omega_{X_1/k}^1\tensor \ccO_X\lra \Omega_X^1\lra 0.
\end{equation}
Let $f:\P^1_k\ra X$ be a non-constant morphism. We can pull back the exact sequence
(\ref{ses-i}) to $\P^1_k$ and get 
\begin{equation}\label{ses-ii}
0\lra f^* \ccO_X\stackrel{\dt d\eps}{\lra} f^*(\Omega_{X_1/k}^1\tensor \ccO_X)\lra 
          f^* \Omega_X^1\lra 0
\end{equation}
\begin{claim}
If (\ref{ses-ii}) is a non-split sequence, then $f$ does not extend to $X_1$.
\end{claim}

\begin{pf}
Suppose $f$ extends to an $A$-linear morphism  $f_1:\P^1_A\ra X_1$. Then we get a 
commutative diagram 
\begin{eqnarray}
\begin{array}{ccccccc}
0\lra & f^* \ccO_X  &\stackrel{\dt d\eps}{\lra}& f^*(\Omega_{X_1/k}^1\tensor \ccO_X) 
                          &\lra& f^* \Omega_X^1 &\lra 0\\
      &\alpha\Big\da&    &       \beta\Big\da    &    & \Big\da          &      \\
0\lra&\ccO_{\P^1_k}&\stackrel{\dt d\eps}{\lra}& \Omega_{\P^1_A/k}^1\tensor\ccO_{\P^1_k}   
                          &\lra&\Omega_{\P^1_k}^1&\ra 0
\end{array}
\end{eqnarray}
Since $f_1$ is an $A$-linear morphism, it follows that $\beta(d\eps)=d\eps$. This
implies that $\alpha$ is non-zero and therefore $\alpha$ is an isomorphism. 
Since $\P^1_A$ is a trivial extension of $\P^1_k$, the bottom row is a split exact 
sequence. Therefore the top row splits as well. 
\end{pf} 

\begin{pf}[Proof of (\ref{fo-lemma-ii})]
In view of the last claim it is enough to prove that the short exact sequence 
(\ref{ses-ii}) is non-split for the deformation $X_1=X_{\ccL}$ of $X$ if $\ccL$ is
ample. For that it is enough to show that the image of 
$\iota_{\sigma}^{-1}[\ccL]\in \H^1(X,\ccT_X)$ under the natural map  
$$\H^1(X,\ccT_X)\lra \H^1(\P^1_k,f^*\ccT_X)$$
is non-zero. After applying the isomorphism $\iota_{\sigma}:\ccT_X\ra \Omega_X^1$
it is enough to show that the image of $[\ccL]\in \H^1(X,\Omega_X^1)$ under the
natural map 
$$\H^1(X,\Omega_X^1)\lra \H^1(\P^1_k,f^*\Omega_X^1)$$
is non-zero. But the image of $[\ccL]$ under the composition 
$$\H^1(X,\Omega_X^1)\lra \H^1(\P^1_k,f^*\Omega_X^1)\lra \H^1(\P^1,\Omega_{\P^1_k}^1)$$
is just $[f^*\ccL]\in \H^1(\P^1,\Omega_{\P^1}^1)$. Since $f^*\ccL$ is ample,
it follows that $[f^*\ccL]$ is non-zero. This finishes the proof of the proposition.
\end{pf}

\section{Proof of Theorem \ref{gencoiso-thm}}

\subsection{A Cohomological Lemma}

Let $f:V\ra W$ be a morphism of separated schemes of finite type over
$k$. For each $i\ge 0$ we can construct a homomorphism
$$f^*:\H^p(W,\ccO_W)\lra \H^p(V,\ccO_V)$$
as follows. Let $\fW=\{W_i\}_i$ be an open affine cover of $W$ and let 
$\fV=\{V_i\}_i$ be an open affine cover of $V$, that is a refinement of 
$f^{-1}\fW$. Then we define $f^*$ to be the composition 
$$\H^p(W,\ccO_W)\cong \cH^p(\fW,\ccO_W)\ra \cH^p(f^{-1}\fW,\ccO_V)
  \ra \cH^p(\fV, \ccO_V)\cong \H^i(V,\ccO_V). $$

\begin{lem}\label{funct-lem}
Let $f:V\ra W$ be a projective morphism of smooth, projective and complex 
varieties $V$ and $W$. Then the diagram
\begin{eqnarray*}
\begin{CD}
\H^p(W,\ccO_W) @>{\cong}>> \overline{\H^0(W,\Omega_W^p)} \\
@V{f^*}VV @V{f^*_{p,0}}VV \\
\H^p(V,\ccO_V) @>{\cong}>> \overline{\H^0(V,\Omega_V^p)}
\end{CD}
\end{eqnarray*}
commutes. Here $f^*$ is the map defined above and $f^*_{p,0}$ is the pull back
of an $(p,0)$-form. The horizontal maps are the the isomorphisms given by 
Hodge theory and $ ^-$ denotes complex conjugation.
\end{lem}

\begin{pf}
Since the isomorphisms  
$H^{0,p}_{\overline{\del}}=\ccH^{0,p}=\overline{\ccH^{p,0}}
=\overline{\H^{p,0}}$
obviously commute with pull back of forms, it is enough to 
show that the diagram 
\begin{eqnarray*}
\begin{CD}
\H^p(W,\ccO_W) @>{\cong}>> \H^{0,p}_{\overline{\del}}(W) \\
@V{f^*}VV @V{f^*_{0,p}}VV \\
\H^p(V,\ccO_V) @>{\cong}>> \H^{0,p}_{\overline{\del}}(W)
\end{CD}
\end{eqnarray*}
commutes, where the horizontal maps are the Dolbeault isomorphisms.

Let $\fW=\{W_i\}_i$ be an open affine cover of $W$. For any sheaf $\ccF$ on $W$
we introduce the notation
$$\Gamma(\fW^a,\ccF):=\prod_{W_{i_1},\dots,W_{i_a}\in \fW} 
  \Gamma(W_{i_1}\cap \cdots \cap W_{i_a},\ccF).$$
Note that $\Gamma(\fW^{\bullet},\ccF)$, together with the appropriate 
differential maps, is a {\v C}ech resolution of $\ccF$.
Let $\ccC_W^{\bullet}$ be the complex of sheaves on $W$ with $\ccO_W$ in degree
$0$ and zero in all other degrees. Let $\ccA_W^{0,p}$ be the sheaf of 
$(0,p)$-forms on $W$. Then $(\ccA_W^{0,\bullet},\overline{\del})$ is
a complex of sheaves and the morphism of sheaves 
$\overline{\del}:\ccO_W\ra \ccA_W^{0,1}$ induces a morphism of complexes
$\ccC^{\bullet}_W\ra \ccA_W^{0,\bullet}.$
By applying $\Gamma(\fW^{\bullet},-)$ we obtain a morphism of double complexes
$$\{\Gamma(\fW^a, \ccC^b_W)\}_{a,b} \lra\{\Gamma(\fW^a, \ccA^{0,b}_W)\}_{a,b}.$$
Let $\fV$ be an open affine cover of $V$, that is a refinement of $f^{-1}\fW$
and do the same construction we just did with $W$ now for $V$.
Since by [GH, p.24],
$\overline{\del}_V \circ f^*= f^*\circ \overline{\del}_W$
we obtain a commutative diagram of double complexes
\begin{eqnarray*}
\begin{CD}
\{\Gamma(\fW^a, \ccC^b_W)\}_{a,b} @>>> \{\Gamma(\fW^a, \ccA^{0,b}_W)\}_{a,b} \\
@VVV @VVV \\
\{\Gamma(\fV^a, \ccC^b_V)\}_{a,b} @>>> \{\Gamma(\fV^a, \ccA^{0,b}_V)\}_{a,b}
\end{CD}
\end{eqnarray*}
Now we can take the cohomology of the total complexes and obtain a commutative
diagram
\begin{eqnarray*}
\begin{CD}
\H^p(W, \ccO_W) @>>> \H^{0,p}_{\overline{\del}}(W) \\
@V{f^*}VV @V{f^*_{0,p}}VV \\
\H^p(V, \ccO_V) @>>> \H^{0,p}_{\overline{\del}}(V) \\
\end{CD}
\end{eqnarray*}
where the horizontal maps are the Dolbeault-isomorphisms. This finishes the 
proof of the lemma.
\end{pf}

\subsection{Proof of Theorem \ref{gencoiso-thm}}

Let $X,Y,f,D,D_i,E_i,r,\sigma$ as in Theorem \ref{gencoiso-thm}.
We prove the theorem in steps. For any $2$-form $\alpha$ on a variety we define
its {\em rank} to be the biggest number $2j$, such that $\alpha^{\wedge j}$ is 
not generically zero. 

\begin{claim}\label{claim-linalg}
$\rank(\sigma|D_i)\ge 2d_i-2n$.
\end{claim}

\begin{pf}
Let $\eta$ be the generic point of $D_i$. At $\eta$ we interpret 
$\sigma$ as a non-degenerate alternating bilinear form on 
the $k(\eta)$-vector space $\ccT_X\tensor k(\eta)$. 
Then $\sigma|D_i$ is at $\eta$ obtained as the 
restriction of the bilinear form on the subspace
$\Theta_D\tensor k(\eta)\subset \ccT_X\tensor k(\eta)$. It is now easy to see
that $\rank (\sigma|D_i)\ge 2d_i-2n$.
\end{pf}

\begin{claim}\label{claim-movcurv}
$2d_i-2n\ge e_i$.
\end{claim}

\begin{pf}
By [Kaw, Theorem 1], we know that the exceptional fibres of $f:X\ra Y$ 
are uniruled. 
Let $F_i\subset X$ be a generic fibre of $D_i\ra F_i$. 
Let $\phi:\P^1\ra F_i$ be a rational curve that cannot be bent and broken
in $F_i$. Therefore 
$$\dim_{[\phi]}Mor(\P^1,F_i)\le 2 \dim F_i+1=2d_i-2e_i+1.$$
Since the image of $\phi$ in $X$ gets contracted under $f:X\ra Y$, all its 
deformations in $X$ stay in the exceptional set and we may assume that
all small deformations stay in $D_i$. Therefore 
\begin{eqnarray*}
\dim_{[\phi]} Mor(\P^1,X)&=&\dim_{[\phi]} Mor(\P^1,D_i) \\
&=&\dim_{[\phi]} Mor(\P^1,F_i)+\dim E_i\le 2d_i-e_i+1.
\end{eqnarray*}
Since $X$ is symplectic, it follows from (\ref{cor-ann-mov}), that 
$\dim_{[\phi]}Mor(\P^1,X)\ge 2n+1.$
Therefore we obtain that $2d_i-e_i+1\ge 2n+1$ and hence
$2d_i-2n \ge e_i.$ 
\end{pf}

Note that Claim \ref{claim-linalg} and Claim \ref{claim-movcurv}
imply Theorem \ref{gencoiso-thm}.(i). 
Let us now suppose that $X$ is projective. Theorem \ref{gencoiso-thm}.(ii)
is an immediate consequence of the next claim.

\begin{claim}\label{claim-hodge}
$\rank (\sigma |D_i)\le e_i$.
\end{claim}

\begin{pf}
Let $D_i'\ra D_i$ be a resolution. The commutative diagram 
\begin{eqnarray*}
\begin{array}{ccc}
D_i'&\lra & X \\
\da & & \da \\
E_i &\lra & Y
\end{array}
\end{eqnarray*}
induces a commutative diagram
\begin{eqnarray*}
\begin{array}{ccc}
\H^p(D_i',\ccO_{D_i'})& \longleftarrow & \H^p(X,\ccO_X) \\
\uparrow & & \uparrow \\
\H^p(E_i,\ccO_{E_i}) & \longleftarrow & \H^p(Y,\ccO_Y)
\end{array}
\end{eqnarray*}
Since $\omega_X\cong \ccO_X$, Grauert-Riemenschneider vanishing implies
$R^pf_*\ccO_X=0$ for all $p\ge 1$. Therefore the Leray spectral sequence 
with respect to $f:X\ra Y$ degenerates and we obtain that the map 
$f^*:\H^p(Y,\ccO_Y)\ra \H^p(X,\ccO_X)$
is an isomorphism. Since $\H^p(E_i,\ccO_{E_i})=0$ for $p>e_i$, we obtain that
$$\H^p(X,\ccO_X) \lra \H^p(D_i',\ccO_{D_i'})$$
is the zero map for $p>e_i$. By (\ref{funct-lem}), it follows that 
$$\H^0(X,\Omega_X^p) \lra \H^0(D_i',\Omega_{D_i'}^p)$$
is the zero map for $p>e_i$. In particular, $(\sigma|D_i')^j$ is zero for
$2j>e_i$. This shows that $\rank (\sigma|D_i)\le e_i$.
\end{pf}

Therefore we must have $2d_i-2n=e_i$. It also follows that each rational curve 
on $F_i$ moves in a family of dimension at least $2\dim F_i+1$. Then 
[CMSB, Theorem 1.1] implies that the normalization of 
$F_i$ consists of a number of projective spaces $\P^{d_i-e_i}$.
This finishes the proof of Theorem \ref{gencoiso-thm}.

\section{Divisorial Contractions I}\label{div-cont-i}

Let $(X,\sigma)$ be a quasi-projective symplectic manifold of dimension $2n$
and let $f:X\ra Y$ be a projective birational contraction. Let $D\subset X$
and $E\subset Y$ be the exceptional sets with their reduced structure. 
We suppose that $E$ is irreducible and $(2n-2)$-dimensional. This implies
that $D$ is a Cartier divisor. Let $D=D_1\cup\cdots\cup D_r$ be the 
decomposition into irreducible components. We suppose further that $f$ has only
$1$-dimensional exceptional fibres.

\begin{prop}\label{tree-of-lines} 
(i) All exceptional fibres of $f$ with their reduced structure are trees
of $\P^1$'s. \\
(ii) A generic scheme theoretic exceptional fibre of $f$ is an 
A-D-E configuration of $\P^1$'s.
\end{prop}

\begin{pf}
(i) Let $F$ be an exceptional fibre with its reduced structure. Since 
$\omega_X\cong \ccO_X$, Grauert-Riemenschneider vanishing implies 
that $R^1 f_* \ccO_X=0$ and therefore $\H^1(F,\ccO_F)=0$. This forces $F$ to
be a tree of $\P^1$'s. \\
(ii) We cut down $Y$ with $(2n-2)$ generic hyperplane sections to obtain 
a surface $Y'$. The inverse image $X'=f^{-1}(Y')$ is smooth. 
Since $f$ is crepant it follows that $X'\ra Y'$ is a crepant resolution. 
Therefore it is a resolution of a Du-Val singularity. 
Its exceptional set is an A-D-E configuration of $\P^1$'s.
\end{pf}

\begin{prop}
Let $F\subset X$ be a generic exceptional fibre.
Any two irreducible components of $F\cap D_i$ are deformation equivalent.
\end{prop}

\begin{pf}
Since $F$ is a tree of $\P^1$'s, we can find an irreducible component 
$A\subset F\cap D_i$ that intersects at most one other irreducible component 
of $F\cap D_i$. We may suppose that $F$ is an exceptional set of a surface
resolution $S\ra T$ as in the proof of the last proposition. Then 
$$(A.D_i)_X=(A.F\cap D_i)_S$$
and we conclude that $(A.D_i)=-2$ if $A=F\cap D_i$ and $(A.D_i)=-1$ otherwise.
Let $B$ be a small deformation of $A$ (i.e. such that $A\cong B$). 
Then $(B.D_i)<0$ and therefore $B\subset D_i$, i.e. all small deformations
of $A$ stay in $D_i$. Since $A$ moves in a $(2n-2)$-dimensional family in $X$
it follows that the deformations of $A$ dominate $D_i$. In particular,
the deformations $B$, such that $B\cong A$, fill out an open, dense subset
of $D_i$. 
Since $F$ was assumed to be generic, we may assume that $F\cap D_i$ is 
contained in this subset. This shows that $A$ is deformation equivalent
to any other irreducible component of $F\cap D_i$.
\end{pf}

\begin{prop}\label{classif-gen}
Let $F\subset X$ be a generic exceptional fibre with dual graph $\Gamma$.
Then $\Gamma$ admits a (possibly trivial) permutation group $G$ of 
automorphisms, such that two irreducible components of $F$ are in the same 
irreducible component of $D$ if and only if their vertices in $\Gamma$ lie 
in the same $G$-orbit. In particular, there is a 1-1 correspondence
 $$\{\text{irred. components of $D$}\} \longleftrightarrow
   \{\text{$G$-orbits of $\Gamma$}\}.$$ 
The cases $(\Gamma, G)$ are as follows: $(A_l,\{\id\})$, $(D_l,\{\id\})$, 
$(E_l,\{\id\})$, $(A_l,S_2)$, $(D_l,S_2)$, $(D_4,S_3)$, $(E_6,S_2)$. 
\end{prop}

\begin{pf}
$\Gamma$ is of type $A_l, D_l$ or $E_l$.
We make a case by case analysis for the types of $\Gamma$. 
Suppose that $\Gamma$ is of type $A_l$. We label a vertex
with $i$, if the curve corresponding to the vertex is contained in $D_i$.
The case $(A_l,\{\id\})$ corresponds to the labelling 
$$1-2-3-\cdots-l.$$
Now suppose that we have started off labelling with 
$$1-2-3-\cdots-(k-1)-k-\bullet-*-\cdots$$
In view of the last proposition we conclude that $\bullet$ can only be 
labelled with $k+1$, $k$, or $k-1$. The first case 
gives nothing new. In the two other cases it is easy to see that we have
to continue as follows:
$$1-2-3-\cdots-(k-1)-k-k-(k-1)-(k-2)-(k-3)-\cdots$$
$$1-2-3-\cdots-(k-1)-k-(k-1)-(k-2)-(k-3)-(k-4)-\cdots$$
Note that the graph must end on the right side just when we reach the label 
$1$ again. Summarizing, we see that $A_{2l}$ can be labelled 
$1-\cdots-2l$ or $1-\cdots-l-l-\cdots-1$ and $A_{2l-1}$ can be labelled
$1\cdots-(2l-1)$ or $1-\cdots-l-\cdots-1$. 
We conclude that for $\Gamma=A_l$ we only obtain the cases $(A_l,\{\id\})$ and
$(A_l,S_2)$ as stated in the theorem. The generator of $S_2$ acts on $A_l$
by the obvious involution. This finishes the proof of the theorem
for the case $\Gamma=A_l$. The other cases are dealt with using similar 
arguments. We leave the details to the reader.
\end{pf}

\begin{prop}
Let $F\subset X$ be a generic exceptional fibre and let $C_i$ be some 
irreducible component of $F\cap D_i$. Then
$N_1(X/Y)=\oplus_{i=1}^r \R[C_i]$.
\end{prop}

\begin{pf}
$N_1$ is generated by $f$-exceptional irreducible curves. 
Any such curve $B$ is a $\P^1$, which we know must move in a 
$(2n-2)$-dimensional family. Therefore $B$ is deformation equivalent to some
linear combination of components of $F$. Therefore $N_1$ is generated by
the irreducible components of $F$. Since any two components of $F\cap D_i$ are
deformation equivalent, it follows that $N_1$ is generated by the $C_i$'s.
To finish the proof it is enough to show that the intersection matrix
$M=\{(C_i.D_j)_X\}_{ij}$ has non-zero determinant. We do a case by case analysis
with respect to the classification in the last proposition. Denote by
$M(\Gamma, G)$ the intersection matrix, where $\Gamma$ is the dual graph of 
$F$ and $G$ is the permutation group of automorphisms. Then it is easy
to compute that
\begin{eqnarray*}
\begin{array}{lcll}
M(A_l,\{\id\}) &=& C(A_l) & l\ge 1 \\
M(D_l,\{\id\}) &=& C(D_l) & l\ge 4 \\
M(E_l,\{\id\}) &=& C(E_l) & l=6,7,8 \\
M(A_2,S_2)    &=& (1)    &         \\
M(A_{2l},S_2) &=& diag(1,\dots,1,\frac{1}{2}) C(B_l) & l\ge 2 \\
M(A_{2n-1},S_2) &=& C(B_l) & l\ge 2 \\
M(D_l,S_2) &=& C(C_{l-1}) & l\ge 4 \\
M(D_4,S_3) &=& C(G_2) &  \\
M(E_6,S_2) &=& C(F_4) &  
\end{array}
\end{eqnarray*}
where $C(\Delta)$ denotes the Cartan matrix of type $\Delta$ 
(cf. [Hum, p. 59]). Since all Cartan matrices have non-zero determinant it 
follows that all intersection matrices $M$ have non-zero determinant, which 
finishes the proof.
\end{pf}

\begin{prop}\label{cone-thm}
Let $i$ be arbitrary.
(i) For  any sufficiently small $\eps>0$, $(X,\eps D_i)$ is a klt-pair. \\
(ii) $[C_i]$ is an $f$-relative extremal ray for $(X,\eps D_i)$. \\
(iii) There exists a factorization $X\stackrel{g_i}{\ra} Z_i\ra Y$, such
that $D_i$ is the $g_i$-exceptional set.
\end{prop}

\begin{pf}
(i) Let $\phi:\hat{X}\ra X$ be an embedded resolution of $D_i$ 
with strict transform $Q\subset \hat{X}$. 
Let $M_j\subset \hat{X}$ be the exceptional divisors. We have 
$$K_{\hat{X}}\equiv \phi^*K_{X}+\sum a_j M_j$$
$$Q \equiv \phi^* D_i+\sum b_j M_j$$
for some $a_j$ and $b_j$. Note that all $a_j$ are $>0$, since $X$ is smooth.
Combining both equations we obtain 
$$K_{\hat{X}}+\eps Q \equiv \phi^*(K_{X}+\eps D_i)+\sum (a_j+\eps b_j) M_j.$$
Recall that $(X,\eps D_i)$ is klt if 
$a_j+\eps b_j>-1, \quad \text{all $i$}$
and $\lfloor \eps Q \rfloor \le 0$. Both conditions are satisfied for $\eps>0$
sufficiently small. \\
(ii) This is an immediate consequence from the last proposition.\\
(iii) This follow from (i), (ii) and the relative cone theorem [KM, p. 95].
\end{pf}

\section{Irreducible Divisorial Contractions}

\subsection{Basic Facts}\label{irr-basic-facts}

We keep $f:X\ra Y, D,E,n$ as in the last section and assume furthermore
that $D$ is irreducible. Let $F\subset X$ be a generic exceptional fibre
with dual graph $\Gamma$. Then (\ref{classif-gen}) implies that $\Gamma$ is 
either $A_1$ or $A_2$.

\begin{prop}\label{normal-bundle}
Let $C\cong \P^1$ be an irreducible component of some exceptional fibre. Then 
$$\ccN_{C/X}\cong \ccO(-2)\oplus \ccO(-1)^{a} 
      \oplus \ccO^{2n-2a-2} \oplus \ccO(1)^{a}.$$ 
\end{prop}

\begin{pf}
Let $\ccI\subset \ccO_X$ be the ideal sheaf defining $C\subset X$. 
Take cohomology in the exact sequence
$$0\lra \ccN_{C/X}^{\vee}\lra \ccO_X/\ccI^2\lra \ccO_C\lra 0$$
to obtain 
$$\H^0(\ccO_X/\ccI^2)\twoheadrightarrow \H^0(\ccO_C)
\lra \H^1(\ccN_{C/X}^{\vee})\lra \H^1(\ccO_X/\ccI^2).$$
Since by Grauert-Riemenschneider vanishing, $R^1f_*\ccO_X=0$, it follows
that $\H^1(\ccO_X/\ccI^2)=0$ and therefore $\H^1(\ccN_{C/X}^{\vee})=0$.
Now take cohomology in the exact sequence
\begin{eqnarray}
0\lra \ccN_{C/X}^{\vee}\lra \Omega_X|C\lra \Omega_C\lra 0 \label{n-b-seq}
\end{eqnarray}
to obtain that $h^1(\Omega_X|C)=1$. Let 
$\Omega_X|C\cong \sum\nolimits_{i=1}^{2n}\ccO_C(\alpha_i)$
be the decomposition with $\alpha_1\le \cdots \le \alpha_{2n}$. 
Then $h^1(\Omega_X|C)=1$ implies that $\alpha_1=-2$ and $\alpha_2\ge -1$. 
But since $X$ is symplectic we have $\Omega_X|C\cong (\Omega_X|C)^{\vee}$ 
and therefore $\alpha_i=\alpha_{2n-i}$ for all $i$.
Therefore $\Omega_X|C$ must be of the form
$$\ccO(-2)\oplus\ccO(-1)^a\oplus\ccO^{2n-2a-2}\oplus\ccO(1)^a\oplus\ccO(2).$$
Now it follows from the sequence (\ref{n-b-seq})
that $\ccN_{C/X}$ is of the desired form.
\end{pf}

Let $H\subset \Hilb(X)$ be the irreducible component of the Hilbert scheme that
contains an irreducible component of some generic $f$-exceptional fibre.
Let $\pi:U\ra H$ be the universal family. We have the following diagram

\begin{eqnarray*}
\begin{array}{ccccc}
U & \stackrel{\phi}{\lra} & D & \subset & X \\
\Big\downarrow\vcenter{\rlap{$\scriptstyle \pi$}} & & 
  \Big\downarrow\vcenter{\rlap{$\scriptstyle g$}} & & 
  \Big\downarrow\vcenter{\rlap{$\scriptstyle f$}} \\
H & \stackrel{\psi}{\lra} & E & \subset & Y.
\end{array}
\end{eqnarray*}

\begin{prop}\label{typei-big-prop-ii}
(i) All scheme theoretic fibres of $\pi:U\ra H$ are reduced $\P^1$'s.\\
(ii) $\pi:U\ra H$ is a smooth morphism.\\
(iii) $H$ is smooth and $2n-2$-dimensional.\\
(iv) $U$ is smooth.
\end{prop}

\begin{pf}
(i) Any irreducible, reduced exceptional curve $C$ is a $\P^1$ and it moves in
$X$ in a family of dimension at least $2n-2$. Therefore the unbroken 
deformations of any such $C$ dominate $D$. Let $B$ be some deformation of 
$C$ and $B'\subset B$ an irreducible reduced component. The unbroken 
deformations of $B'$ dominate $D$ and therefore $B'$ is deformation equivalent
to $C$. Since $C$ is also deformation equivalent to $B$, it follows that $B'$ 
is deformation equivalent to $B$. This is only possible if $B'=B$. This shows 
that all deformations of $C$ are reduced $\P^1$'s. \\
(ii) Since $\pi:U\ra H$ is flat and all fibres are smooth, it follows
that $\pi$ is smooth.\\
(iii) By (\ref{mov-thm}), we know that each rational curve in $X$ 
moves in an at least $(2n-2)$-dimensional family on $X$. Therefore 
$\dim H\ge 2$. But (\ref{normal-bundle}) implies that for each $[C]\in H$ the 
tangent space $T_{H_1,[C]}=H^0(C, \ccN_{C/X})$
is $(2n-2)$-dimensional. Therefore $H$ is smooth and $(2n-2)$-dimensional.\\
(iv) By (iii) $H$ is smooth and by (ii) $\pi$ is smooth. Therefore $U$ is 
smooth.
\end{pf}

\begin{lem}\label{typei-phi-finite}
The morphism $\phi:U\ra D$ is finite. 
\end{lem}

\begin{pf}
Suppose not. Then there is a reduced, irreducible curve $A\subset U$, such that
$\phi(A)$ is a point. Since all reduced fibres of $\pi:U\ra H$ are mapped 
isomorphically via $\phi$ to $D$, it follows that $B=\pi(A)$ must be a curve. 
All $\pi$-fibres over $B$ have images under $\phi$ that pass through the point 
$p=\phi(A)$. Therefore they are all contained in the $f$-exceptional fibre 
$F=f^{-1}f(p)$, i.e. we have a map
$$\phi|_{\pi^{-1}(B)}: \pi^{-1}(B)\lra F.$$
Since $\dim F=1$ and $\dim \pi^{-1}(B)=2$ and for each $b\in B$, the map
$$\phi|_{\pi^{-1}(b)}: \pi^{-1}(b)\lra \phi(\pi^{-1}(b))$$
is an isomorphism, we can find two points $b,b'\in B$, such that 
$\phi(\pi^{-1}(b))=\phi(\pi^{-1}(b')).$
But this contradicts the universal property of the Hilbert scheme.
\end{pf}

Let $Y'$ be obtained from $Y$ by taking $2n-3$ generic hyper plane sections
and let $U',H', X', D',E'$ be obtained by the base change $Y'\ra Y$.
Since the hyper planes are generic we may assume the following.
\begin{itemize}
\item $Y'$ is an irreducible and normal $3$-fold.
\item $X'$ is a smooth irreducible $3$-fold.
\item $H'$ is a smooth curve (possibly with several connected components).
\item $U'$ is a smooth surface (possibly with several connected components).
\item $E'$ is an irreducible and reduced curve.
\end{itemize}
Since $f:X\ra Y$ is a crepant contraction, adjunction implies that $X'\ra Y'$
is also a crepant.

\subsection{The Case $(\Gamma,G)=(A_1,\{\id\})$}

Let us assume that $(\Gamma,G)=(A_1,\{\id\})$.

\begin{thm}\label{a-i-nondeg}
$D$ and $E$ are smooth. $f|_D:D\ra E$ is a smooth morphism all whose fibres 
are $\P^1$'s.
\end{thm}

Before we give the proof we want to recall the following result of Wilson.

\begin{thm}\label{wilson}
Let $V\ra V'$ be a crepant contraction of a smooth, projective $3$-fold, that
contracts an irreducible surface $S\subset V$ to a curve $C$. 
Then $C$ is a smooth curve.
\end{thm}

\begin{pf} 
The result is stated and proved in the case of a contraction
from a projective Calabi-Yau threefold in [Wi, 3.1]. 
Note, however, that the proof is 
entirely local, so the assumption that $V$ be projective and Calabi-Yau
can be weakened to the assumption that the contraction be crepant.
\end{pf}

\begin{pf}[Proof of Theorem \ref{a-i-nondeg}.]
Since $\Gamma=A_1$ it follows that a generic exceptional fibre is a single
$\P^1$. Therefore $\psi:H\ra E$ is birational. Since $E'$ is irreducible
this implies that $H'$ is irreducible and therefore also $U'$ and $D'$.
In particular, the birational contraction $f':X'\ra Y'$ satisfies the 
conditions of (\ref{wilson}). We deduce that $E'$ is smooth. Since $E'$ is 
a generic curve on $E$, this implies that $E$ is regular in codimension $1$. 
Let $E^{sm}\subset E$ be the regular locus and $D^{sm}\subset D$ its inverse 
image in $D$. Since $D^{sm}$ is 
Cohen-Macaulay, $E^{sm}$ smooth and $f|_{D^{sm}}:D^{sm}\ra E^{sm}$ 
equidimensional, it follows that $f|_{D^{sm}}:D^{sm}\ra E^{sm}$
is flat. Therefore it exhibits a family of curves in  $X$, which, by the 
universal property of the Hilbert scheme, yields a map $D^{sm}\ra U$,
such that the composition $D^{sm}\ra U\stackrel{\phi}{\ra} D$ is the embedding.
Since $U$ is smooth, this implies that $D^{sm}$ is smooth and therefore $D$ 
is regular in codimension $1$. Since $D$ is Cohen-Macaulay it follows that 
$D$ must be normal. But $\phi:U\ra D$ is the normalization map. Therefore
$U\cong D$ and in particular $D$ is smooth. Consider the short exact sequence
$$0\lra \ccO_X(-D)\lra \ccO_X\lra \ccO_D\lra 0.$$
Since $f$ is crepant and $\ccO_X(-D)$ is $f$-ample, vanishing implies 
$R^1f_* \ccO_X(-D)=0$. Therefore the map $f_*\ccO_X \ra f_*\ccO_D$ is 
surjective. But $f_*\ccO_X=\ccO_Y$ and $f_*\ccO_D=\ccO_{\tilde{E}}$, where 
$\tilde{E}\ra E$ is the normalization. Therefore the natural map 
$\ccO_Y\ra \ccO_{\tilde{E}}$ is surjective. This implies that $E$ is in fact 
normal. Since $\psi:H\ra E$ is the normalization,
it follows that $H\cong E$ and that $E$ is smooth. In particular we have
$(U\ra H)\cong (D\ra E)$ and therefore $D\ra E$ is a smooth morphism all whose
fibres are reduced $\P^1$'s.
\end{pf}

\begin{cor}\label{omega-e-triv-a-i}
$\omega_{E}\cong \ccO_{E}$.
\end{cor}

\begin{pf}
Let $\ccI=\ccO_X(-D)$. Since $\ccI$ is $f$-ample,
it follows that in the diagram
\begin{eqnarray*}
\begin{array}{ccccccccc}
0 & \lra & \ccI/\ccI^2 & \lra & \Omega_X|D & \lra & \Omega_{D}           & \lra & 0 \\
  &     &             &     & ||         &     &                      &         \\
0 & \lra & \ccT_{D}    & \lra & \ccT_X|D   & \lra & (\ccI/\ccI^2)^{\vee} & \lra & 0
\end{array}
\end{eqnarray*}
the composition $\ccI/\ccI^2  \ra (\ccI/\ccI^2)^{\vee}$ is zero and 
that we get a surjective map $\Omega_{D} \ra (\ccI/\ccI^2)^{\vee}.$
The composition 
$$f^*\Omega_{E}\lra \Omega_{D} \lra (\ccI/\ccI^2)^{\vee}$$
is zero, since for any fibre $C\cong \P^1$ of $D\ra E$, we have
$f^*\Omega_{E}|C\cong \ccO_C^{\oplus 2n-2}$ and 
$(\ccI/\ccI^2)^{\vee}|C\cong \ccO_C(-2)$. 
Therefore we obtain an isomorphism 
$\Omega_{D/E} \cong (\ccI/\ccI^2)^{\vee}.$
We can now calculate as follows
\begin{eqnarray*}
f^*\det \Omega_{E} 
  & \cong & \det\Omega_{D}\tensor (\det \Omega_{D/E})^{\vee} \\
  & \cong & \omega_{D} \tensor (\ccI/\ccI^2) \\
  & \cong & \ccO_X(D)|D\tensor \ccO_X(-D)|D \\
  & \cong & \ccO_{D}.
\end{eqnarray*}
Finally the projection formula 
$f_* f^*\omega_{E}\cong f_*\ccO_{D}\tensor \omega_{E}$
implies that $\omega_{E}\cong \ccO_{E}$.
\end{pf}

\subsection{The Case $(\Gamma,G)=(A_2,S_2)$}

Let us now assume that $(\Gamma,G)=(A_2,S_2)$. We aim to prove the 
following result.

\begin{thm}\label{a-ii-nondeg}
After possibly deleting a closed subscheme of codimension $\ge 4$ in $Y$
the following is true.
$E$ is smooth, $\psi:H\ra E$ is \'etale of degree $2$. $\pi:U\ra H$ has
a section $\Sigma \subset U$, such that $\Sigma\ra \phi(\Sigma)$ is \'etale
of degree $2$ and $f:\phi(\Sigma)\ra E$ is an isomorphism. $\Sigma$ is 
precisely the locus in $U$, where $\phi:U\ra D$ is not an isomorphism.
$D\ra E$ is flat and all fibres are $A_2$-configurations of $\P^1$'s.
\end{thm}

\begin{cor}\label{omega-e-triv-a-ii}
$\omega_E\cong \ccO_E$.
\end{cor}

\begin{pf}[Proof of Corollary.]
Consider the diagram 
\begin{eqnarray*}
\begin{array}{ccccccc}
0\lra& \ccT_{\Sigma} & \lra & \ccT_U|\Sigma & \lra & \ccN_{\Sigma/U} & \lra 0 \\
     &               &      &      ||       &      &                 &       \\
0\lra& \ccT_{U/H}|\Sigma &\lra& \ccT_U|\Sigma &\lra& \pi^*\ccT_H|\Sigma & \lra 0.
\end{array}
\end{eqnarray*}
Since $\Sigma\cong H$, it follows that the composition 
$\ccT_{\Sigma}\ra\ccT_U|\Sigma\ra  \pi^*\ccT_H|\Sigma$
is an isomorphism. Therefore the composition
$\ccT_{U/H}|\Sigma \ra  \ccT_U|\Sigma \ra  \ccN_{\Sigma/U}$
is an isomorphism. It follows that we have the following commutative diagram
\begin{eqnarray*}
\begin{array}{ccc}
\ccT_{U/H} & \lra & \ccT_{U/H}|\Sigma \\
\da        &     &  ||               \\
\ccT_U     & \lra & \ccN_{\Sigma/U}
\end{array}
\end{eqnarray*}
which can be completed to the following diagram 
\begin{eqnarray}\label{tgt-diag}
\begin{array}{ccccccc}
     & 0                   &    & 0          &    &                   &       \\
     &  \da                &    & \da        &    &                   &       \\
0\lra& \ccT_{U/H}(-\Sigma) &\lra& \ccT_{U/H} &\lra& \ccT_{U/H}|\Sigma & \lra 0 \\
     &  \da                &    & \da        &    &  ||                       \\
0\lra& \ccT_U(-\log\Sigma) &\lra& \ccT_U     &\lra& \ccN_{\Sigma/U}   & \lra 0 \\
     &  \da                &    & \da        &    &                   &       \\
     & \pi^*\ccT_H         &=   &\pi^*\ccT_H &    &                   &       \\
     &  \da                &    & \da        &    &                   &       \\
     & 0                   &    & 0.         &    &                   &       
\end{array}
\end{eqnarray}
The symplectic $2$-form 
$\sigma$ on $X$ gives an isomorphism $\Omega_X\cong \ccT_X$,
from which we can get the following diagram 
\begin{eqnarray*}
\begin{array}{ccccccc}
0\lra & \ccN_{D/X}^{\vee} & \lra & \Omega_X|D & \lra & \Omega_D & \lra 0 \\
     &            &     & ||         &     &          &       \\
0\lra & \theta_D   & \lra & \ccT_X|D   & \lra & \ccN_{D/X} & \lra 0
\end{array}
\end{eqnarray*}
Since $\ccO_X(-D)$ is $f$-ample, it follows that the composition
$\ccN_{D/X}^{\vee} \ra \ccN_{D/X}$ 
must be zero. Therefore we obtain a map $\ccN_{D/X}^{\vee}\ra \theta_D$. 
Pulling this back to $U$ yields a non-zero map
$\phi^* \ccN_{D/X}^{\vee}\ra \phi^*\theta_D\ra \ccT_U(-\log \Sigma).$
Let $C\subset U$ be a $\pi$-fibre. Then 
$\phi^* \ccN_{D/X}^{\vee}|C\cong \ccO_C(1)$.
Since $\pi^*\ccT_H|C\cong \ccO^{2n-2}_C$ it follows that the composition, 
obtained by the diagram (\ref{tgt-diag}), 
$\phi^* \ccN_{D/X}^{\vee}\ra\ccT_U(-\log \Sigma)\ra \pi^*\ccT_H$
is zero and therefore we have a non-zero map
$$\phi^* \ccN_{D/X}^{\vee}\lra\ccT_{U/H}(-\Sigma),$$
which is an isomorphism, since $\ccT_{U/H}(-\Sigma)|C\cong \ccO_C(1)$. 
Since $\omega_X\cong \ccO_X$, it follows from adjunction that 
$\omega_D=\ccN_{D/X}$. Since $D$ is obtained from $U$ by glueing along 
$\Sigma$, it follows that $\omega_U=\phi^*\omega_D\tensor \ccO_U(-\Sigma)$. 
Therefore
\begin{eqnarray*}
\pi^*\omega_H &=& \omega_U\tensor \omega_{U/H}^{\vee} \\
              &=& \phi^* \omega_D \tensor \ccO_U(-\Sigma)\tensor \omega_{U/H}^{\vee} \\
              &=& \phi^* \ccN_{D/X} \tensor \ccO_U(-\Sigma)\tensor \ccT_{U/H}
\end{eqnarray*}
But we already know that $\phi^* \ccN_{D/X}^{\vee}\cong\ccT_{U/H}(-\Sigma)$
and therefore $\pi^*\omega_H\cong \ccO_U$. Pushing down to $H$ gives
$\omega_H\cong \ccO_H$.

Since $\phi(\Sigma)\cong E$ and $H\cong \Sigma$, we see that 
$\tau:=\sigma^{\wedge n-1}|\phi(\Sigma)$ gives a non-zero global section of
$\omega_E$. Since $\omega_H$ is trivial and $H\ra E$ \'etale, this implies 
that $\omega_E$ is trivial.
\end{pf}

Let us now turn to the proof of the theorem.
Since $\phi:U\ra D$ exhibits the normalization of $D$, we can consider the
conductor ideal $\fc\subset \ccO_{U}$ with respect to $\phi$. It defines a 
closed subscheme $\Xi\subset U$, which is of pure codimension $1$, since 
$D$ satisfies Serre's condition $S_2$. Since the generic fibre of $D\ra E$ 
consists of two $\P^1$'s that meet in a single reduced point, it follows that 
we can find an irreducible component $\Sigma\subset \Xi$, such that 
$\pi|_{\Sigma}:\Sigma\ra H$ is dominant. Since the meeting point of the two 
$\P^1$'s is reduced, it follows that $\Sigma$ is generically reduced. Since 
$\Sigma$ is of pure codimension $1$ in $U$ it follows that it is reduced.

The morphism $\Sigma\ra H$ is proper and birational. 
Let $Y^0\subset Y$ be an open subset and 
$\Sigma^0, \Xi^0, U^0, H^0, X^0, D^0, E^0$ be obtained by the base change
$Y^0\ra Y$. We choose $Y^0$ maximal, such that $\Sigma^0\ra H^0$ is 
an isomorphism. Since $H$ is smooth, it follows that the locus in $H$, where
$\Sigma\ra H$ is not an isomorphism is of codimension at least $2$. In
particular, the complement of $Y^0\subset Y$ has codimension at least $4$.

Since $E'\subset E$ is a generic curve on $E$, we may assume that 
$E\subset E^0$.

Let $T=U^0\oplus^{\Sigma^0} \phi(\Sigma^0)$ be obtained from $U^0$ by glueing 
$\Sigma^0\subset U^0$ with $\Sigma^0\ra \phi(\Sigma^0)$. By definition, we 
can factor $\phi:U^0\ra D^0$ into 
$$U^0\stackrel{\tau}{\lra} T \lra D^0.$$

\begin{claim}\label{claim-t-d-iso}
$T\ra D^0$ is an isomorphism. In particular, $\Sigma=\Xi$.
\end{claim}

\begin{pf}
Since $D^0$ satisfies $S_2$, it is enough to show that $T\ra D^0$ is an 
isomorphism in codimension $1$. Let $T', \Sigma'$ be obtained by the base
change $Y'\ra Y$. It is now enough to show that $T'\ra D'$ is an isomorphism.
Note that $T'$ is obtained from $U'$ by glueing $\Sigma'\ra \phi(\Sigma')$.

(i) Suppose $D'$ is irreducible. Then we can apply (\ref{wilson}) to
$f':X'\ra Y'$. We deduce
that $E'$ is smooth and all fibres of $D'\ra E'$ are two $\P^1$'s that meet 
in a point. In particular, the normalization $\tilde{D'}\ra D'$ is a smooth 
ruled surface $\tilde{D'}\ra C$, where 
$C\ra E'$ is an \'etale cover of degree $2$. Since $\tilde{D'}\ra C$ exhibits 
a family of rational curves in $X$ there is a map $\tilde{D'}/C\ra U/H$. 
Since $C\ra H'$ is birational and $H'$ is smooth, it follows that $C\ra H'$ is 
an isomorphism. Therefore we can obtain $D'$ from $U'$ by glueing
$\Sigma'\ra \phi(\Sigma')$. This shows that $T'\ra D'$ is an isomorphism.

(ii) Suppose $D'$ is reducible. Then $D'=D'_1\cap D'_2$ with $D'_i$ 
irreducible. Let $W\ra \overline{W}_i$ be the morphism that contracts $D'_i$ 
to a curve, say $C_i$. We can apply (\ref{wilson}) and
deduce that $C_i$ is smooth. 
Therefore $D'_i\ra C_i$ is flat with generic fibre a single $\P^1$. 
Since these $\P^1$'s cannot degenerate it follows that all fibres are $\P^1$'s 
and  that $D'_i\ra C_i$ is a smooth ruled surface. In particular, 
the definition of the Hilbert scheme implies that there is a morphism 
$D'_i/C_i\ra U/H$. Since $H'$ is smooth and $H'\ra E'$ of degree 
$2$ it follows that $C_1\sqcup C_2\ra H'$ is an isomorphism. 
Since $C_i\ra E'$ is birational, it follows that $C_1\cong C_2$. 
Let $H'=H'_1\sqcup H'_2$ with $C_i\ra H'_i$. Let $U'_i=\pi^{-1}(H'_i)$. 
Let $\Sigma'_i=U'_i\cap \Sigma'$. Let $S$ be 
obtained by glueing $U'_1$ and $U'_2$ along $\Sigma'_1\cong \Sigma'_2$. 
Then we have maps 
$$U'\lra S \lra \tau(S)\lra D'.$$
Note that the composition $S\ra F$ is an isomorphism in codimension $1$. 
Since $D'$ satisfies $S_2$ it follows that $S\ra D'$ is an isomorphism. 
Since $S\cong T'$ by construction, it follows that $T\ra D'$ is an isomorphism.
\end{pf}

Let $\Delta=\phi(\Sigma)$ and $\Delta^0=\Delta\times_Y Y^0$.

\begin{claim}
$\Delta^0$ is normal.
\end{claim}

\begin{pf}
Suppose not. Since $D^0=U^0\oplus^{\Sigma^0}\Delta^0$ we get a morphism
$\alpha:U^0\oplus^{\Sigma^0}\tilde{\Delta}^0\ra D^0$,
where $\tilde{\Delta}^0$ is the normalization of $\Delta^0$, which is an 
isomorphism in codimension $1$. Since $D^0$ satisfies $S_2$ it follows
that $\alpha$ must be an isomorphism and therefore 
$\tilde{\Delta}^0\cong \Delta^0$.
\end{pf}

\begin{claim}
$\Sigma^0\ra \Delta^0$ is \'etale of degree $2$ over the smooth locus of 
$\Delta^0$.
\end{claim}

\begin{pf}
The proof of (\ref{claim-t-d-iso}) implies that $\Sigma^0\ra \Delta^0$ is 
\'etale of degree $2$ in codimension $1$. Therefore the purity of the branch 
locus implies that $\Sigma^0\ra \Delta^0$ is \'etale over the smooth locus of
$\Delta^0$.
\end{pf}

\begin{claim}
$\Delta^0$ is smooth.
\end{claim}

\begin{pf}
Suppose not. We may assume that the singular locus of $\Delta^0$ contains a 
component of dimension $k$. Let $\overline{Y}\subset Y$ be obtained by taking 
$k$ generic hyper plane sections and let $\overline{X}, \overline{D}, 
\overline{E}$, etc. be obtained by the base change $\overline{Y}\ra Y$.
Then $\overline{\Sigma}^0\ra \overline{\Delta}^0$ is \'etale of degree
$2$ outside a finite number of points. Let $p\in \overline{\Delta}^0$ 
be a singular point. Then locally analytically at $p$, 
$\overline{\Sigma}^0\ra \overline{\Delta}^0$ can be described by 
$$k[[x_1,\dots,x_l]]^{\Z/2\Z}\hookrightarrow k[[x_1,\dots,x_l]],
\quad l=\dim \overline{\Delta}^0,$$
where the generator of $\Z/2\Z$ acts via $x_i\mapsto -x_i$. Since
$k[[x_1,\dots,x_l]]^{\Z/2\Z}=k[[x_i x_j | i,j=1,\dots,n]]$ it follows
that $\overline{\Delta}^0$ has embedding dimension $\frac{1}{2} l(l+1)$. 
But $l=2n-2-k$ and $\overline{\Delta}^0\subset \overline{X}^0$ with 
$\overline{X}^0$ smooth, it follows that 
$$\frac{1}{2} l(l+1)\le \dim \overline{X}^0=l+2.$$
This is only satisfied for $l\le 2$. Since $\overline{\Delta}^0$ is normal 
it follows that $l\ge 2$. Therefore $l=2$. 

Let $h\in H$ be the point, such that $p\in\phi(\pi^{-1}(h))$. Let
$A=\ccO_{H,h}$. Since $\overline{U}^0\ra \overline{H}^0$
is a $\P^1$-bundle with a section $\overline{\Sigma}^0$, the bundle is locally
trivial in the Zariski topology and therefore 
$\pi^{-1}(\Spec A)=\P^1\times \Spec A$. We may assume that $\overline{U}^0$ is
on an open affine patch over $\Spec A$ given by $\Spec A[z]$ and 
$\overline{\Sigma}^0$ is defined by $(z=0)$. Let $B=A^{\Z/2\Z}$ be the ring
of invariants of the $\Z/2\Z$-action, i.e. $B=\ccO_{\overline{\Delta}^0,\psi(h)}$.

We can describe the cofibred sum square 
\begin{eqnarray*}
\begin{array}{ccc}
\overline{U}^0 & \lra & \overline{D}^0 \\
\cup & \oplus & \cup \\
\overline{\Sigma}^0 & \lra & \overline{\Delta}^0 
\end{array}
\end{eqnarray*}
in a neighbourhood of $p$ by  
\begin{eqnarray*}
\begin{array}{ccc}
A[y] & \longleftarrow  & R \\
\big\downarrow\vcenter{\rlap{$\scriptstyle{y=0}$}}  &       &\big\downarrow \\
A   & \longleftarrow  & B
\end{array}
\end{eqnarray*}
In other words,
\begin{eqnarray*}
R&=&\{ \sum a_i z^i \, | \, \sum a_i z^i \mod z \in B \} \\
  & =&\{ b+fz \, | \, b\in B, f\in A[z]\}.
\end{eqnarray*}
Let $a_1,a_2\in A$ generate $\m_A/\m_A^2$. 
Let $b_1, b_2,b_3\in B$ generate $\m_B/\m_B^2$. Then we have 
$$z, a_1 z, a_2 z, b_1,b_2,b_3\in R.$$
Clearly, these elements are independent generators of $\m_R/\m_R^2$.
Therefore the embedding dimension of $\overline{D}^0$ is at least six.
This is a contradiction.
\end{pf}

\begin{claim}
$E^0\cong \Delta^0$. In particular, $E^0$ is smooth.
\end{claim}

\begin{pf}
Since $U^0, \Sigma^0$ and $\Delta^0$ are smooth, it follows that 
$D^0=U^0\oplus^{\Sigma^0}\Delta^0\ra E^0$ factors through the normalization
$\tilde{E}^0$ of $E^0$. In particular, $D^0\ra \tilde{E}^0\ra E^0$
is the Stein factorization of $D^0\ra E^0$. Since $\Delta^0\ra E$ is birational
and finite and since $\Delta^0$ is smooth, it is enough to prove that $E^0$ is normal. For that it is enough to show that 
$$\ccO_{Y^0}=f_*\ccO_{X^0}\lra f_*\ccO_{D^0}\cong \ccO_{\tilde{E}^0}$$
is surjective. For this it is enough to prove that $R^1f_*\ccO_{X^0}(-D^0)=0$.
But $\ccO_{X}(-D)$ is $f$-ample and $f$ is crepant, therefore vanishing
implies $R^1f_*\ccO_{X^0}(-D^0)=0$.
\end{pf}

\section{Divisorial Contractions II}

We consider a contraction $f:X\ra Y$ as in section \ref{div-cont-i}.
We assume the notations and results of that section. We have classified
such $f$ in terms of $(\Gamma, G)$, where $\Gamma$ is an $A-D-E$ diagram
and $G$ is a permutation group of automorphisms of $\Gamma$.
If $f$ is of type $(A_{2l}, S_2)$ we allow to delete a closed subscheme of 
codimension $\ge 4$ in $Y$.

\begin{prop}
(i) The normalization $\tilde{E}$ of $E$ is smooth. \\
(ii) For any $i,j$,  $D_i\cap D_j$ is either empty or smooth of dimension
$2n-2$ and $D_i\cap D_j\ra \tilde{E}$ is \'etale. \\
(iii) $D\ra E$ factors through $\tilde{E}$. \\
(iv) $D\ra \tilde{E}$ is a flat morphism with constant fibres.
\end{prop}

\begin{pf}
By (\ref{cone-thm}), we can for each component $D_i$ of $D$ find a 
factorization $X\ra Z_i\ra Y$, such that $X\ra Z_i$ contracts only $D_i$.
The contraction $X\ra Z_i$ is irreducible and either of type $(A_1,\{\id\})$
or type $(A_2,S_2)$. Therefore we can apply (\ref{a-i-nondeg}) or
(\ref{a-ii-nondeg}) respectively to obtain that the locus of singularities
$F_i\subset Z_i$ is smooth and that $D_i\ra F_i$ is a flat bundle all whose
fibres are either $\P^1$'s or two copies of $\P^1$'s that meet in a single 
point. Moreover (\ref{omega-e-triv-a-i}) or (\ref{omega-e-triv-a-ii}) implies
that $\omega_{F_i}$ is trivial. In the following we will have to do a case
by case analysis with respect to the $(\Gamma,G)$-classification.

(i) An inspection of all cases reveals that we can always find some component
of $D$, say $D_1$, such that $F_1\ra E$ is birational. Since this morphism
is also finite and $F_1$ is smooth, it follows that $F_1\ra E$ exhibits the 
normalization of $E$ and in particular, $\tilde{E}$ is smooth.

(ii) Suppose $S=D_i\cap D_j$ is nonempty. Since $D_1$ and $D_2$ are
divisors in $X$ it follows that $S$ is Cohen-Macaulay and $(2n-2)$-dimensional.
Since the $f$-exceptional curves of $D_i$ and $D_j$ are not numerically 
equivalent, it follows that $S$ does not contain $f$-exceptional curves.
Therefore $S\ra E$ is a finite and surjective morphism. By looking at a 
generic $f$-exceptional fibre we conclude that $S$ is generically reduced. 
Since it is Cohen-Macaulay it is reduced. An inspection of all
cases reveals that at least for one of $i$ and $j$, say $i$ we have that
$S\ra F_i$ is birational. Since $F_i$ is smooth it follows that $S\ra F_i$
is an isomorphism and in particular $S$ is smooth and $\omega_S$ is trivial.
Since $S$ is smooth it follows that $S\ra E$ factors through $\tilde{E}$. 
Since $\omega_{\tilde{E}}$ is also trivial it follows that $S\ra \tilde{E}$
is \'etale.

(iii) $D$ is obtained by glueing the various $D_i$ along $D_i\cap D_j$. Since
all morphisms $D_i\ra E$ and $D_i\cap D_j\ra E$ factor through $\tilde{E}$
it follows that $D\ra E$ factors through $\tilde{E}$.

(iv) Since $D$ is Cohen-Macaulay, $D\ra \tilde{E}$ is equidimensional and
$\tilde{E}$ is smooth, it follows that $D\ra \tilde{E}$ is flat. 
By (\ref{a-i-nondeg}) and (\ref{a-ii-nondeg}) it follows that $D_i\ra F_i$ 
have constant fibres and since $F_i\ra \tilde{E}$ is \'etale, it follows
that $D_i\ra \tilde{E}$ has constant fibres. Therefore all fibres of
$D\ra \tilde{E}$ have the same number of components. In order to show that
the fibres are constant it remains to show that the intersection behaviour
of the fibre components stays constant. But since any fibre component 
$C\subset D_i$ can be deformed to a component $C'$ of a generic fibre it 
follows that $(C.D_j)_X=(C'.D_j)_X$. Therefore the intersection behaviour
stays constant.
\end{pf}

\begin{prop}
$E$ is normal.
\end{prop}

\begin{pf}
Since $D\ra E$ factors through $\tilde{E}$ it is enough to show that 
$$\ccO_Y\cong f_*\ccO_X\lra f_*\ccO_D\cong \ccO_{\tilde{E}}$$
is surjective. Therefore it is enough to show that $R^1f_*\ccO_X(-D)=0$.
Recall that, since $f$ is crepant, we have $R^1f_*\ccO(A)=0$ for any 
$f$-nef $A$.
We do a case by case analysis with respect to $(\Gamma,G)$.

$(A_l,\{\id\})$. It is easy to check that $-D$ is $f$-nef. Therefore
$R^1f_*\ccO_X(-D)=0$ and we are done.

$(A_{2l},S_2)$. It is easy to check that $-D$ is $f$-nef. Therefore
$R^1f_*\ccO_X(-D)=0$ and we are done.

$(D_l,\{\id\})$. We label the components of $D$ as\\ 
\smallskip \\
\centerline{
\setlength{\unitlength}{0.0005in}
{\em
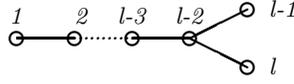
\begin{picture}(2708,776)(1118,-754)
\thicklines
\put(1201,-361){\circle{150}}
\put(1801,-361){\circle{150}}
\put(2401,-361){\circle{150}}
\put(3001,-361){\circle{150}}
\put(3601,-61){\circle{150}}
\put(3601,-661){\circle{150}}
\put(1201,-361){\line( 1, 0){600}}
\put(2401,-361){\line( 1, 0){600}}
\put(3601,-61){\line(-2,-1){600}}
\put(3001,-361){\line( 2,-1){600}}
\multiput(1801,-361)(60.0,0.0){11}{\makebox(6.6667,10.0){.}}
\put(1126,-211){\makebox(0,0)[lb]{\smash{1}}}
\put(1801,-211){\makebox(0,0)[lb]{\smash{2}}}
\put(2251,-211){\makebox(0,0)[lb]{\smash{l-3}}}
\put(2851,-211){\makebox(0,0)[lb]{\smash{l-2}}}
\put(3826,-136){\makebox(0,0)[lb]{\smash{l-1}}}
\put(3826,-736){\makebox(0,0)[lb]{\smash{l}}}
\end{picture}
}}
\smallskip \\
Let $A_j=-D-\sum_{i=j+1}^{l-2}D_i$ for $j=1,\dots,l-2$. It is easy to check 
that $A_1$ is $f$-nef, therefore $R^1f_*\ccO_X(A_1)=0$. Consider
$$0\ra \ccO_X(A_{j-1})\ra \ccO_X(A_j)\ra\ccO_{D_j}(A_j\cap D_j)\ra 0.$$
It is easy to see that $\ccO_{D_j}(A_j\cap D_j)$ has degree $\ge -1$ on 
each $f$-exceptional curve in $D_j$. Therefore 
$R^1f_*\ccO_{D_j}(A_j\cap D_j)=0$ and induction on $j$ implies that 
$R^1f_*\ccO_X(A_j)=0$ for all $j$. Since $A_{l-2}=-D$ we are done.

$(D_l,S_2)$. For this and the remaining other cases we apply the same 
procedure as above. We specify divisors $A_1,\dots,A_m$, such that $A_1$
is $f$-nef, $A_m=-D$ and such that for all $j$ we have 
$A_{j+1}-A_j=D_{i_j}$ for some $i_j$. Moreover we will also have that 
$\ccO_{D_{i_j}}(A_{j+1}\cap D_{i_j})$ has degree $\ge -1$ on each 
$f$-exceptional curve. This allows us to conclude that $R^1f_*\ccO_X(-D)=0$.
We label the components of $D$ as follows.\\
\smallskip\\
\centerline{
\setlength{\unitlength}{0.0005in}
{\em
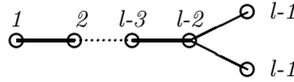
\begin{picture}(2708,785)(1118,-763)
\thicklines
\put(1201,-361){\circle{150}}
\put(1801,-361){\circle{150}}
\put(2401,-361){\circle{150}}
\put(3001,-361){\circle{150}}
\put(3601,-61){\circle{150}}
\put(3601,-661){\circle{150}}
\put(1201,-361){\line( 1, 0){600}}
\put(2401,-361){\line( 1, 0){600}}
\put(3601,-61){\line(-2,-1){600}}
\put(3001,-361){\line( 2,-1){600}}
\multiput(1801,-361)(60.00000,0.00000){11}{\makebox(6.6667,10.0000){.}}
\put(1126,-211){\makebox(0,0)[lb]{\smash{1}}}
\put(1801,-211){\makebox(0,0)[lb]{\smash{2}}}
\put(2251,-211){\makebox(0,0)[lb]{\smash{l-3}}}
\put(2851,-211){\makebox(0,0)[lb]{\smash{l-2}}}
\put(3826,-136){\makebox(0,0)[lb]{\smash{l-1}}}
\put(3826,-736){\makebox(0,0)[lb]{\smash{l-1}}}
\end{picture}
}}
\smallskip\\
We define the $A_j$'s as $A_j=-D-\sum_{i=j+1}^{l-2}D_i$ for $j=1,\dots l-2$.

$(D_4,S_3)$. We label the components of $D$ as follows.\\
\smallskip\\
\centerline{
\setlength{\unitlength}{0.0005in}
{\em 
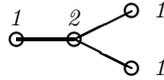
\begin{picture}(1508,785)(2318,-763)
\thicklines
\put(2401,-361){\circle{150}}
\put(3001,-361){\circle{150}}
\put(3601,-61){\circle{150}}
\put(3601,-661){\circle{150}}
\put(2401,-361){\line( 1, 0){600}}
\put(3601,-61){\line(-2,-1){600}}
\put(3001,-361){\line( 2,-1){600}}
\put(2326,-211){\makebox(0,0)[lb]{\smash{1}}}
\put(2926,-211){\makebox(0,0)[lb]{\smash{2}}}
\put(3826,-136){\makebox(0,0)[lb]{\smash{1}}}
\put(3826,-736){\makebox(0,0)[lb]{\smash{1}}}
\end{picture}
}}
\smallskip\\
We define $A_1=-D_1-2 D_2$, $A_2=-D$.

$(E_6,\{\id\})$. We label the components of $D$ as follows.\\
\smallskip\\
\centerline{
\setlength{\unitlength}{0.0005in}
{\em
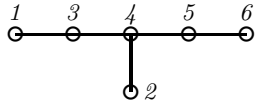
\begin{picture}(2566,960)(2318,-1063)
\thicklines
\put(2401,-361){\circle{150}}
\put(3001,-361){\circle{150}}
\put(4201,-361){\circle{150}}
\put(4801,-361){\circle{150}}
\put(3601,-961){\circle{150}}
\put(3601,-361){\circle{150}}
\put(2401,-361){\line( 1, 0){2400}}
\put(3601,-361){\line( 0,-1){600}}
\put(2326,-211){\makebox(0,0)[lb]{\smash{1}}}
\put(3730,-1036){\makebox(0,0)[lb]{\smash{2}}}
\put(2926,-211){\makebox(0,0)[lb]{\smash{3}}}
\put(3526,-211){\makebox(0,0)[lb]{\smash{4}}}
\put(4126,-211){\makebox(0,0)[lb]{\smash{5}}}
\put(4726,-211){\makebox(0,0)[lb]{\smash{6}}}
\end{picture}
}}
\smallskip\\
We define $A_1=-D_1-2D_2-2D_3-3D_4-2D_5-D_6$, $A_2=A_1+D_2$, $A_3=A_2+D_4$,
$A_4=A_3+D_3$, $A_5=A_4+D_5$, $A_6=A_5+D_4$. Note that $A_6=-D$.

$(E_6,S_2)$. We label the components of $D$ as follows. \\
\smallskip\\
\centerline{
\setlength{\unitlength}{0.0005in}
{\em
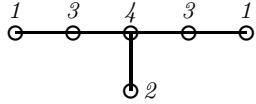
\begin{picture}(2566,960)(2318,-1063)
\thicklines
\put(2401,-361){\circle{150}}
\put(3001,-361){\circle{150}}
\put(4201,-361){\circle{150}}
\put(4801,-361){\circle{150}}
\put(3601,-961){\circle{150}}
\put(3601,-361){\circle{150}}
\put(2401,-361){\line( 1, 0){2400}}
\put(3601,-361){\line( 0,-1){600}}
\put(2326,-211){\makebox(0,0)[lb]{\smash{1}}}
\put(3730,-1036){\makebox(0,0)[lb]{\smash{2}}}
\put(2926,-211){\makebox(0,0)[lb]{\smash{3}}}
\put(3526,-211){\makebox(0,0)[lb]{\smash{4}}}
\put(4126,-211){\makebox(0,0)[lb]{\smash{3}}}
\put(4726,-211){\makebox(0,0)[lb]{\smash{1}}}
\end{picture}
}}
\smallskip\\
We define $A_1=-D_1-2D_2-2D_3-3D_4$, $A_2=A_1+D_2$, $A_3=A_2+D_4$, 
$A_4=A_3+D_3$, $A_5=A_4+D_4$. Note that $A_5=-D$.

$(E_7,\{\id\})$. We label the components of $D$ as follows. \\
\smallskip\\
\centerline{
\setlength{\unitlength}{0.0005in}
{\em
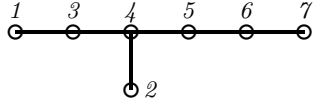
\begin{picture}(3166,960)(2318,-1063)
\thicklines
\put(2401,-361){\circle{150}}
\put(3001,-361){\circle{150}}
\put(4201,-361){\circle{150}}
\put(4801,-361){\circle{150}}
\put(3601,-961){\circle{150}}
\put(3601,-361){\circle{150}}
\put(5401,-361){\circle{150}}
\put(2401,-361){\line( 1, 0){2400}}
\put(3601,-361){\line( 0,-1){600}}
\put(4801,-361){\line( 1, 0){600}}
\put(2326,-211){\makebox(0,0)[lb]{\smash{1}}}
\put(3730,-1036){\makebox(0,0)[lb]{\smash{2}}}
\put(2926,-211){\makebox(0,0)[lb]{\smash{3}}}
\put(3526,-211){\makebox(0,0)[lb]{\smash{4}}}
\put(4126,-211){\makebox(0,0)[lb]{\smash{5}}}
\put(4726,-211){\makebox(0,0)[lb]{\smash{6}}}
\put(5326,-211){\makebox(0,0)[lb]{\smash{7}}}
\end{picture}
}}
\smallskip\\
We define $A_1=-D_1-2D_2-3D_4-4D_5-3D_6-2D_7$. To define the other $A_j$'s
we add successively the components
$$D_7,D_6,D_5,D_2,D_4,D_5,D_3,D_4,D_6,D_5.$$

$(E_8,\{\id\})$. We label the components of $D$ as follows. \\
\smallskip\\
\centerline{
\setlength{\unitlength}{0.0005in}
{\em
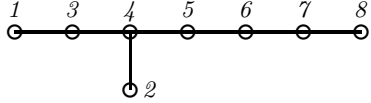
\begin{picture}(3766,960)(2318,-1063)
\thicklines
\put(2401,-361){\circle{150}}
\put(3001,-361){\circle{150}}
\put(4201,-361){\circle{150}}
\put(4801,-361){\circle{150}}
\put(3601,-961){\circle{150}}
\put(3601,-361){\circle{150}}
\put(5401,-361){\circle{150}}
\put(6001,-361){\circle{150}}
\put(2401,-361){\line( 1, 0){2400}}
\put(3601,-361){\line( 0,-1){600}}
\put(4801,-361){\line( 1, 0){600}}
\put(5401,-361){\line( 1, 0){600}}
\put(2326,-211){\makebox(0,0)[lb]{\smash{1}}}
\put(3730,-1036){\makebox(0,0)[lb]{\smash{2}}}
\put(2926,-211){\makebox(0,0)[lb]{\smash{3}}}
\put(3526,-211){\makebox(0,0)[lb]{\smash{4}}}
\put(4126,-211){\makebox(0,0)[lb]{\smash{5}}}
\put(4726,-211){\makebox(0,0)[lb]{\smash{6}}}
\put(5326,-211){\makebox(0,0)[lb]{\smash{7}}}
\put(5926,-211){\makebox(0,0)[lb]{\smash{8}}}
\end{picture}
}}
\smallskip\\
We define $A_1=-2D_1-3D_2-4D_3-6D_4-5D_5-4D_6-3D_7-2D_8$. To define
the other $A_j$'s we add successively the components
\begin{eqnarray*}
\begin{array}{l}
D_8,D_7,D_6,D_5,D_4,D_2,D_3,D_4,D_1,D_5,D_3,D_6,D_4,D_7,D_5,D_2, \\
D_6,D_4,D_5,D_3,D_4.
\end{array}
\end{eqnarray*}
We have now dealt with all possible cases. This finishes the proof.
\end{pf}

\section{Global Structure of the Singular Locus}

Let $(X,\sigma)$ be a projective, symplectic manifold of dimension $2n$ and
let $f:X\ra Y$ be a projective, birational contraction. Suppose $f$ contracts
a divisor. Then by (\ref{gencoiso-thm}), we find an irreducible component
$E$ of the singular locus of $Y$, such that $\dim E=2n-2$. 
Let $D$ be an irreducible component of $f^{-1}(E)_{red}$ 
and let $\tilde{E}$ be the normalization of $E$.

\begin{thm}
(i) $\omega_{\tilde{E}}\cong \ccO_{\tilde{E}}$.\\
(ii) $\tilde{E}$ is a symplectic variety.
\end{thm}

\begin{pf}
(i) Let $Y_1\subset Y$ be obtained by deleting all points $y\in Y$ with 
$\dim f^{-1}(y)\ge 2$. In case of type $(A_{2l},\Z_2)$ we allow to delete
a closed subscheme of $Y_1$ of codimension $\ge 4$. 
Let $f_1, X_1, E_1$ be obtained by the base change $Y_1\ra Y$. 
Then $f_1:X_1\ra Y_1$ satisfies in a neighbourhood of $E_1$ the
conditions of (\ref{big-thm}) and therefore we can conclude that $E_1$ is 
smooth and $\omega_{E_1}$ is trivial. Since the complement of $E_1\subset E$
has codimension at least $2$ it follows that $\omega_{\tilde{E}}$ is also
trivial.

(ii) Let $E'\ra \tilde{E}$ be a resolution of singularities and let 
$D'\ra D\times_E E'$ a resolution of singularities. Let $g:D'\ra E'$ and
$h:D'\ra X$. We get a commutative diagram
\begin{eqnarray*}
\begin{array}{ccc}
\H^2(D',\ccO_{D'}) & \longleftarrow & \H^2(X,\ccO_X) \\
\big\uparrow           &            & \big|\big| \\
\H^2(E',\ccO_{E'}) & \longleftarrow & \H^2(Y,\ccO_Y).
\end{array}
\end{eqnarray*}
We can apply (\ref{funct-lem}) to obtain that for some 
$\eps\in \H^0(E',\Omega^2_{E'})$ we have $g^*\eps=h^*\sigma$. Since $\sigma$
is symplectic, it follows that $h^* \sigma^{n-1}$ is non-zero. Therefore
$\eps^{n-1}$ is non-zero. Since $\omega_{\tilde{E}}$ is trivial, $\eps$
is non-degenerate on the smooth locus of $\tilde{E}$. Therefore $\tilde{E}$
is a symplectic variety.
\end{pf}

\section{Examples}

\begin{ex}
The following is well known. See [Slo] or [Hin] for a reference.
Let $G$ be a simple algebraic group and $B\subset G$ a Borel subgroup.
Let $X=T^*(G/B)$ and let $Y$ be the nilpotent cone of $\mathfrak g$, where 
$\mathfrak g$ is the lie algebra of $G$. There is a birational morphism
$f:X\ra Y$ and $Sing(Y)$ has codimension $2$ in $Y$. Let $f_1:X_1\ra Y_1$
be the restriction that deletes all fibres of dimension $\ge 2$. Since
$X$ is a cotangent bundle, it carries a natural symplectic form. Therefore
$f_1$ satisfies the conditions of (\ref{big-thm}) and we can associate a type
$(\Gamma, S_m)$ to it. The type of $G$ determines the type of $f_1$. This
is known as Brieskorn's theorem and the correspondence is as follows.
\begin{eqnarray*}
\begin{array}{cc}
\text{Type of $G$} & \text{Type of $f_1$} \\
A_l                & (A_l,\{\id\})         \\
B_l                & (A_{2l-1},S_2)      \\
C_l                & (D_{l+1},S_2)       \\
D_l                & (D_l,\{\id\})         \\
E_l                & (E_l,\{\id\})         \\
F_4                & (E_6,S_2)           \\
G_2                & (D_4,S_3)         
\end{array}
\end{eqnarray*}
Note that this gives examples for all types of contractions except 
$(A_{2l},S_2)$.
\end{ex}

\begin{ex}
We are going to construct an example for a contraction of type $(A_2,S_2)$
from a projective symplectic variety.

Let $E$ be the elliptic curve, which has an automorphism $\eps:E\ra E$ of
order $3$. The automorphism $\eps\times \eps^{-1}$ on $E\times E$ 
has also order $3$ and we let $S=E\times E/\langle\eps\times \eps^{-1}\rangle$.
Let $res:T\ra S$ be a minimal resolution. The origin $(0,0)\in E\times E$ 
is an isolated fixed point of action and its image $s\in S$ is a 
finite quotient singularity of type $\frac{1}{3}(1,2)$, which is a 
Du-Val singularity of type $A_2$. The exceptional fibre $F=res^{-1}(s)$
consists of two $(-2)$-curves meeting in a point. 

Note that $\eps\times \eps^{-1}$ preserves a global generator
of $\omega_{E\times E}$. It follows that $T$ is a $K3$-surface.

Let $f:E\times E\ra E\times E$ be given by $(x,y)\ra (-y,x)$. Then
$f$ is an automorphism of $E\times E$ of order $4$. $f$ preserves 
$\eps\times \eps^{-1}$-orbits and $f$ preserves also a global generator 
of $\omega_{E\times E}$. Therefore $f$ induces automorphisms $g:S\ra S$
and $h:T\ra T$ of order $4$. Note that $h$ preserves a global generator 
of $\omega_T$. Note also that $h$ exchanges the two components of $F$.

Let $A$ be an abelian surface and let $a:A\ra A$ be the translation by
a $4$-torsion point. Then $h\times a$ is an automorphism of order $4$ on 
$T\times A$ and $g\times a$ is an automorphism of order $4$ on $S\times A$.
Starting from $res\times \id:T\times A\ra S\times A$ we can quotient out
the action of $h\times a$ and $g\times a$ to obtain a contraction
$f:X\ra Y$. Since $h\times a$ is fixed point free, it follows that $X$ is 
smooth. Let $\tau$ be a global generator of $\omega_T$ and let 
$\alpha$ be a global generator of $\omega_A$. Let $p$ and $q$ be the 
projections from $T\times A$ to the first and second factor
respectively. Then $p^*\tau+q^*\alpha$ is a symplectic $2$-form on 
$T\times A$, which is invariant under $h\times a$. Therefore $X$ is a 
symplectic $4$-fold. The exceptional fibres of $f$ consist of two $\P^1$'s
that meet in one point. The exceptional set $D$ in $X$ is 
$F\times A /\langle h\times a\rangle$, which is irreducible. Therefore $f:X\ra Y$ is 
a divisorial contraction from a symplectic $4$-fold of type $(A_2,S_2)$.
\end{ex}

\begin{ex}
We are going to construct an example of a contraction of type 
$(A_{2l}, S_2)$. Let $V=\Spec k[x,y]$. Let $\eps:V\ra V$ be the automorphism
given by $x\mapsto\zeta x, y\mapsto \zeta^{-1}y$ for $\zeta$ a primitive
$2l$-th root of unity. The quotient $S=V/\langle\eps\rangle$ is a Du-Val 
singularity of type $A_{2l}$. Let $T\ra S$ be a minimal resolution with 
exceptional set $F\subset T$. 

Let $f:V\ra V$ be the automorphism defined by $x\mapsto y, y\mapsto -x$. 
Then $f$ has order $4$ and preserves the global generator 
$dx\wedge dy\in \Gamma(V,\omega_V)$. It is easy to see that $f$ preserves
$\langle\eps\rangle$-orbits. Therefore $f$ descents to an automorphism $g$ of $S$ and 
can in fact be lifted to an automorphism $h$ of $T$. Note that $h$ preserves 
a global generator of $\omega_T$ and swaps the irreducible components of $F$. 

Let $A$ be an abelian surface and let $a:A\ra A$ be the translation by
a $4$-torsion point. As in the last example, we conclude that 
$$X=T\times A/\langle h\times a\rangle \,\lra\, 
   S\times A/\langle g\times a\rangle=Y$$
is a divisorial contraction of type $(A_{2l},S_2)$ from a symplectic 
manifold $X$.
\end{ex}

\noindent
\vspace{0.5cm}\\
{\bf Acknowledgements.} I am grateful to Nick Shepherd-Barron for his continuous 
support and valuable discussions. I also want to thank Dmitri Kaledin for 
his interest in my work and his valuable comments about Section 2.

\vspace{0.5cm}

\noindent
{\bf Address:} \\
DPMMS/CMS, University of Cambridge, Wilberforce Road, CB3 0WB, United Kingdom. \\
J.Wierzba@dpmms.cam.ac.uk

\end{document}